\documentclass[a4paper, 12 pt]{article}

\usepackage[T2A]{fontenc}
\usepackage[cp1251]{inputenc}
\usepackage[english]{babel}
\usepackage[tbtags]{amsmath}
\usepackage{amsfonts,amssymb}

\usepackage{mathrsfs}

\usepackage{geometry} 
\begin{document}
\title{Riordan-Dirichlet group}
 \author{E. Burlachenko}
 \date{}

 \maketitle
\begin{abstract}
Riordan matrices are infinite lower triangular matrices that correspond to certain operators in the space of formal power series. In this paper, we introduce similar matrices for the space of formal Dirichlet series. We show that these matrices form a group similar to the Riordan group, and we derive  an analog of the Lagrange inversion formula for this group. As an example of the application of these matrices, we obtain an analog of the Abel identities.
\end{abstract}
\section{Introduction}

Riordan matrices (Riordan arrays) are infinite lower triangular matrices that correspond to certain operators in the space of formal power series over the field of real or complex numbers. In Section 3, we introduce similar matrices for the space of formal Dirichlet series. In Section 2, we will note aspects of the theory of Riordan matrices on which we will rely when constructing similar objects. We will call these objects the Riordan-Dirichlet matrices. Analogy between them and the Riordan matrices is so complete that the content of Section 3 almost verbatim repeats the content of  Section 2. But this analogy concerns only the noted aspects. Riordan matrices find wide application in various fields of mathematics while the question of the practical application of the Riordan-Dirichlet matrices is open question. In Section 4, we consider the situation when these matrices prove to be useful for obtaining identities similar to identities obtained by means of Riordan matrices. As an example, we obtain an analog of the Abel identities.
\section{Riordan matrices}

We will associate the columns of matrices with the generating functions of their elements, i.e. with the formal power series. Thus, the expression $Aa\left( x \right)=b\left( x \right)$ means that the column vector multiplied by the matrix $A$ has the generating function $a\left( x \right)$, resultant column vector has the generating function $b\left( x \right)$. We will denoted  $n$th coefficient of the series $a\left( x \right)$ and $n$th row of the matrix $A$ respectively by $\left[ {{x}^{n}} \right]a\left( x \right)$,    $\left[ n,\to  \right]A$.   

Matrix $n$th column of which, $n=0,\text{ }1,\text{ }2,\text{ }...$ ,  has the generating function ${{x}^{n}}a\left( x \right)$ will be  denoted by $\left( a\left( x \right),1 \right)$:

$$\left( a\left( x \right),1 \right)=\left( \begin{matrix}
   {{a}_{0}} & 0 & 0 & 0 & \cdots   \\
   {{a}_{1}} & {{a}_{0}} & 0 & 0 & \cdots   \\
   {{a}_{2}} & {{a}_{1}} & {{a}_{0}} & 0 & \cdots   \\
   {{a}_{3}} & {{a}_{2}} & {{a}_{1}} & {{a}_{0}} & \cdots   \\
   \vdots  & \vdots  & \vdots  & \vdots  & \ddots   \\
\end{matrix} \right).$$
Then
$$\left( a\left( x \right),1 \right)b\left( x \right)=\sum\limits_{n=0}^{\infty }{{{b}_{n}}{{x}^{n}}a\left( x \right)}=\sum\limits_{n=0}^{\infty }{{{x}^{n}}\sum\limits_{m=0}^{n}{{{b}_{m}}{{a}_{n-m}}}}=a\left( x \right)b\left( x \right).$$
Thus, matrix $\left( a\left( x \right),1 \right)$ corresponds to the operator of multiplication by the  series $a\left( x \right)$. Set of all such matrices form an algebra isomorphic to the algebra of formal power series:
$$\left( a\left( x \right),1 \right)+\left( b\left( x \right),1 \right)=\left( a\left( x \right)+b\left( x \right),1 \right),$$
$$\left( a\left( x \right),1 \right)\left( b\left( x \right),1 \right)=\left( b\left( x \right),1 \right)\left( a\left( x \right),1 \right)=\left( a\left( x \right)b\left( x \right),1 \right).$$

Matrix $n$th column of which has the generating function ${{x}^{n}}{{a}^{n}}\left( x \right)$ will be  denoted by $\left( 1,a\left( x \right) \right)$:
$$\left( 1,a\left( x \right) \right)=\left( \begin{matrix}
   1 & 0 & 0 & 0 & 0 & \cdots   \\
   0 & {{a}_{0}} & 0 & 0 & 0 & \ldots   \\
   0 & {{a}_{1}} & a_{0}^{\left( 2 \right)} & 0 & 0 & \cdots   \\
   0 & {{a}_{2}} & a_{1}^{\left( 2 \right)} & a_{0}^{\left( 3 \right)} & 0 & \cdots   \\
   0 & {{a}_{3}} & a_{2}^{\left( 3 \right)} & a_{3}^{\left( 3 \right)} & a_{0}^{\left( 4 \right)} & \cdots   \\
   \vdots  & \vdots  & \vdots  & \vdots  & \vdots  & \ddots   \\
\end{matrix} \right),$$ 
where
$$\left[ {{x}^{n}} \right]{{a}^{m}}\left( x \right)=a_{n}^{\left( m \right)},  \qquad a_{0}^{\left( 0 \right)}=1,  \qquad a_{n}^{\left( 0 \right)}=0,  \qquad a_{n}^{\left( 1 \right)}={{a}_{n}}.$$
Then
$$\left( 1,a\left( x \right) \right)b\left( x \right)=\sum\limits_{n=0}^{\infty }{{{b}_{n}}{{x}^{n}}{{a}^{n}}\left( x \right)=}\sum\limits_{n=0}^{\infty }{{{x}^{n}}\sum\limits_{m=0}^{n}{{{b}_{m}}a_{n-m}^{\left( m \right)}}}=b\left( xa\left( x \right) \right).$$
Thus, matrix $\left( 1,a\left( x \right) \right)$ corresponds to the operator of composition of series. Since
$$\left( 1,a\left( x \right) \right){{x}^{n}}b\left( x \right)={{x}^{n}}{{a}^{n}}\left( x \right)b\left( xa\left( x \right) \right),$$
then
$$\left( 1,a\left( x \right) \right)\left( b\left( x \right),1 \right)=\left( b\left( xa\left( x \right) \right),1 \right)\left( 1,a\left( x \right) \right);$$
since
$$\left( 1,a\left( x \right) \right){{x}^{n}}{{b}^{n}}\left( x \right)={{x}^{n}}{{\left( a\left( x \right)b\left( xa\left( x \right) \right) \right)}^{n}},$$
then
$$\left( 1,a\left( x \right) \right)\left( 1,b\left( x \right) \right)=\left( 1,a\left( x \right)b\left( xa\left( x \right) \right) \right).$$
Denote
$$\left( b\left( x \right),1 \right)\left( 1,a\left( x \right) \right)=\left( b\left( x \right),a\left( x \right) \right).$$
Such matrix is called the Riordan matrix [1] – [6]. $n$th column of the Riordan matrix has the generating function $b\left( x \right){{x}^{n}}{{a}^{n}}\left( x \right)$. If ${{b}_{0}}\ne 0$, ${{a}_{0}}\ne 0$, Riordan matrix is called proper. Proper Riordan matrices form a group, called the Riordan group, whose elements are multiplied by the rule:
$$\left( b\left( x \right),a\left( x \right) \right)\left( f\left( x \right),g\left( x \right) \right)=\left( b\left( x \right)f\left( xa\left( x \right) \right),a\left( x \right)g\left( xa\left( x \right) \right) \right).$$
Identity element of the group is the matrix $\left( 1,1 \right)$. Inverse element for the matrix $\left( b\left( x \right),a\left( x \right) \right)$ is the matrix $\left( {1}/{b\left( x\bar{a}\left( x \right) \right),\bar{a}\left( x \right)}\; \right)$, where
$$\bar{a}\left( x \right)a\left( x\bar{a}\left( x \right) \right)=a\left( x \right)\bar{a}\left( xa\left( x \right) \right)=1.$$
{\bfseries Remark 1.} Alternative notations for the Riordan matrix defining the composition exists. For example, $\left( 1,a\left( x \right) \right)$, ${{a}_{0}}=0$; or $\left( 1,xa\left( x \right) \right)$. We have chosen the notation (see [2] – [4]) similar to notation that appears in Section 3 in the construction of the Riordan-Dirichlet matrices.

In the algebra of formal power series, the power and logarithm of the series $a\left( x \right)$, ${{a}_{0}}=1$, are defined as
$${{a}^{\varphi }}\left( x \right)=\sum\limits_{n=0}^{\infty }{\left( \begin{matrix}
   \varphi   \\
   n  \\
\end{matrix} \right){{\left( a\left( x \right)-1 \right)}^{n}}}, \qquad\log a\left( x \right)=\sum\limits_{n=1}^{\infty }{\frac{{{\left( -1 \right)}^{n-1}}}{n}}{{\left( a\left( x \right)-1 \right)}^{n}}.$$
Power of the series can also be defined as
$${{a}^{\varphi }}\left( x \right)=\exp \left( \varphi \log a\left( x \right) \right)=\sum\limits_{n=0}^{\infty }{\frac{{{\varphi }^{n}}}{n!}}{{\left( \log a\left( x \right) \right)}^{n}}=\sum\limits_{n=0}^{\infty }{\frac{{{s}_{n}}\left( \varphi  \right)}{n!}}{{x}^{n}},$$
where ${{{s}_{n}}\left( \varphi  \right)}/{n!}\;$ are the polynomials in $\varphi $ of degree$\le n$ called convolution polynomials [7]. Explicit form of these polynomials:
$${{s}_{0}}\left( \varphi  \right)=1,  \qquad\frac{{{s}_{n}}\left( \varphi  \right)}{n!}=\sum\limits_{m=1}^{n}{{{\varphi }^{m}}\sum{\frac{b_{1}^{{{m}_{1}}}b_{2}^{{{m}_{2}}}...\text{ }b_{n}^{{{m}_{n}}}}{{{m}_{1}}!{{m}_{2}}!\text{ }...\text{ }{{m}_{n}}!}}}\text{ },$$
where ${{b}_{i}}=\left[ {{x}^{i}} \right]\log a\left( x \right)$ and summation of the coefficient of ${{\varphi }^{m}}$ is over all monomials $b_{1}^{{{m}_{1}}}b_{2}^{{{m}_{2}}}...\text{ }b_{n}^{{{m}_{n}}}$ for which $\sum\nolimits_{i=1}^{n}{i{{m}_{i}}}=n$, $\sum\nolimits_{i=1}^{n}{{{m}_{i}}}=m$.

Matrix of the differential operator in the space of formal power series will be denoted by ${{D}_{x}}$:

$${{D}_{x}}a\left( x \right)={a}'\left( x \right)=\sum\limits_{n=1}^{\infty }{n{{a}_{n}}}{{x}^{n-1}}.$$
Then
$${{\left( a\left( x \right)b\left( x \right) \right)}^{\prime }}=a\left( x \right){b}'\left( x \right)+{a}'\left( x \right)b\left( x \right),  \qquad{{\left( {{a}^{n}}\left( x \right) \right)}^{\prime }}=n{{a}^{n-1}}\left( x \right){a}'\left( x \right),$$
$${{\left( {{a}^{\varphi }}\left( x \right) \right)}^{\prime }}=\varphi {a}'\left( x \right)\sum\limits_{n=1}^{\infty }{\left( \begin{matrix}
   \varphi -1  \\
   n-1  \\
\end{matrix} \right){{\left( a\left( x \right)-1 \right)}^{n-1}}}=\varphi {{a}^{\varphi -1}}\left( x \right){a}'\left( x \right),$$
$${{\left( \log a\left( x \right) \right)}^{\prime }}={a}'\left( x \right)\sum\limits_{n=1}^{\infty }{{{\left( -1 \right)}^{n-1}}{{\left( a\left( x \right)-1 \right)}^{n-1}}}={a}'\left( x \right){{a}^{-1}}\left( x \right).$$
{\bfseries Theorem 1.} \emph{Each formal power series $a\left( x \right)$, ${{a}_{0}}=1$, is associated  by means of  the transform
$${{a}^{\varphi }}\left( x \right)=\sum\limits_{n=0}^{\infty }{\frac{{{x}^{n}}}{{{a}^{\beta n}}\left( x \right)}\left[ {{x}^{n}} \right]}\left( 1-x\beta {{\left( \log a\left( x \right) \right)}^{\prime }} \right){{a}^{\varphi +\beta n}}\left( x \right)\eqno (1)$$
with the set of series $_{\left( \beta  \right)}a\left( x \right)$, $_{\left( 0 \right)}a\left( x \right)=a\left( x \right)$, such that}
$${}_{\left( \beta  \right)}a\left( x{{a}^{-\beta }}\left( x \right) \right)=a\left( x \right),   \qquad a\left( x{}_{\left( \beta  \right)}{{a}^{\beta }}\left( x \right) \right)={}_{\left( \beta  \right)}a\left( x \right),$$
$$\left[ {{x}^{n}} \right]{}_{\left( \beta  \right)}{{a}^{\varphi }}\left( x \right)=\left[ {{x}^{n}} \right]\left( 1-x\beta {{\left( \log a\left( x \right) \right)}^{\prime }} \right){{a}^{\varphi +\beta n}}\left( x \right)=\frac{\varphi }{\varphi +\beta n}\left[ {{x}^{n}} \right]{{a}^{\varphi +\beta n}}\left( x \right),$$
$$\left[ {{x}^{n}} \right]\left( 1+x\beta {{\left( \log {}_{\left( \beta  \right)}a\left( x \right) \right)}^{\prime }} \right){}_{\left( \beta  \right)}{{a}^{\varphi }}\left( x \right)=\frac{\varphi +\beta n}{\varphi }\left[ {{x}^{n}} \right]{}_{\left( \beta  \right)}{{a}^{\varphi }}\left( x \right)=\left[ {{x}^{n}} \right]a_{{}}^{\varphi +\beta n}\left( x \right).$$
{\bfseries Remark 2.} Transformation (1) is a particular case of the Lagrange series expansion. But, considering our plans, we give an alternative proof which relies only on the properties of the Riordan matrices.\\
{\bfseries Proof.}
 If the matrices $\left( 1,{{a}^{-1}}\left( x \right) \right)$, ${{a}_{0}}=1$, $\left( 1,b\left( x \right) \right)$, ${{b}_{0}}=1$, are mutually inverse, then
$$\left( 1,{{a}^{-1}}\left( x \right) \right)b\left( x \right)=a\left( x \right),   \qquad\left( 1,b\left( x \right) \right)a\left( x \right)=b\left( x \right).$$
Since
$${{\left( {{x}^{n}}{{b}^{n}}\left( x \right) \right)}^{\prime }}=n{{x}^{n-1}}{{b}^{n-1}}\left( x \right)b\left( x \right)\left( 1+x{{\left( \log b\left( x \right) \right)}^{\prime }} \right),$$
then
$${{D}_{x}}\left( 1,b\left( x \right) \right)=\left( b\left( x \right)\left( 1+x{{\left( \log b\left( x \right) \right)}^{\prime }} \right),b\left( x \right) \right){{D}_{x}},$$
$$\left( 1,b\left( x \right) \right){a}'\left( x \right)=\frac{{{\left( \log b\left( x \right) \right)}^{\prime }}}{1+x{{\left( \log b\left( x \right) \right)}^{\prime }}}.$$
From this we find:
$${{\left( 1+x{{\left( \log b\left( x \right) \right)}^{\prime }},b\left( x \right) \right)}^{-1}}=\left( 1-x{{\left( \log a\left( x \right) \right)}^{\prime }},{{a}^{-1}}\left( x \right) \right).$$
Denote
$$\left[ {{x}^{n}} \right]{{a}^{m}}\left( x \right)=a_{n}^{\left( m \right)},  \qquad\left[ {{x}^{n}} \right]\left( 1-x{{\left( \log a\left( x \right) \right)}^{\prime }} \right){{a}^{m}}\left( x \right)=c_{n}^{\left( m \right)},$$ 
 $${{a}_{m}}\left( x \right)=\sum\limits_{n=0}^{\infty }{a_{n}^{\left( m+n \right)}}{{x}^{n}},  \qquad{{c}_{m}}\left( x \right)=\sum\limits_{n=0}^{\infty }{c_{n}^{\left( m+n \right)}{{x}^{n}}}.$$
We construct the matrix $A$ $m$th column of which has the generating function ${{x}^{m}}{{a}_{m}}\left( x \right)$ and the matrix $C$ $m$th column of which has the generating function ${{x}^{m}}{{c}_{m}}\left( x \right)$:
$$A=\left( \begin{matrix}
   a_{0}^{\left( 0 \right)} & 0 & 0 & 0 & \cdots   \\
   a_{1}^{\left( 1 \right)} & a_{0}^{\left( 1 \right)} & 0 & 0 & \cdots   \\
   a_{2}^{\left( 2 \right)} & a_{1}^{\left( 2 \right)} & a_{0}^{\left( 2 \right)} & 0 & \cdots   \\
   a_{3}^{\left( 3 \right)} & a_{2}^{\left( 3 \right)} & a_{1}^{\left( 3 \right)} & a_{0}^{\left( 3 \right)} & \cdots   \\
   \vdots  & \vdots  & \vdots  & \vdots  & \ddots   \\
\end{matrix} \right),  \qquad C=\left( \begin{matrix}
   c_{0}^{\left( 0 \right)} & 0 & 0 & 0 & \cdots   \\
   c_{1}^{\left( 1 \right)} & c_{0}^{\left( 1 \right)} & 0 & 0 & \cdots   \\
   c_{2}^{\left( 2 \right)} & c_{1}^{\left( 2 \right)} & c_{0}^{\left( 2 \right)} & 0 & \cdots   \\
   c_{3}^{\left( 3 \right)} & c_{2}^{\left( 3 \right)} & c_{1}^{\left( 3 \right)} & c_{0}^{\left( 3 \right)} & \cdots   \\
   \vdots  & \vdots  & \vdots  & \vdots  & \ddots   \\
\end{matrix} \right).$$
It's obvious that
$$\left[ n,\to  \right]A=\left[ n,\to  \right]\left( {{a}^{n}}\left( x \right),1 \right),$$
$$\left[ n,\to  \right]C=\left[ n,\to  \right]\left( \left( 1-x{{\left( \log a\left( x \right) \right)}^{\prime }} \right){{a}^{n}}\left( x \right),1 \right).$$
Since
$$\left( 1-x{a}'\left( x \right){{a}^{-1}}\left( x \right) \right){{a}^{m}}\left( x \right)={{a}^{m}}\left( x \right)-\frac{x}{m}{{\left( {{a}^{m}}\left( x \right) \right)}^{\prime }},$$
or
$$\left[ {{x}^{n}} \right]\left( 1-x{{\left( \log a\left( x \right) \right)}^{\prime }} \right){{a}^{m}}\left( x \right)=\frac{m-n}{m}\left[ {{x}^{n}} \right]{{a}^{m}}\left( x \right),$$
then
$$\left[ {{x}^{n+m}} \right]A{{x}^{m}}\left( 1-x{{\left( \log a\left( x \right) \right)}^{\prime }} \right){{a}^{-m}}\left( x \right)=\left[ {{x}^{n+m}} \right]C{{x}^{m}}{{a}^{-m}}\left( x \right)=$$
$$=\left[ {{x}^{n}} \right]\left( 1-x{{\left( \log a\left( x \right) \right)}^{\prime }} \right){{a}^{n}}\left( x \right) =1, n=0;\quad =0, n>0.$$
Thus,
$$A=\left( 1+x{{\left( \log b\left( x \right) \right)}^{\prime }},b\left( x \right) \right),   \qquad C=\left( 1,b\left( x \right) \right),$$ 
$$\left[ {{x}^{n}} \right]\left( 1+x{{\left( \log b\left( x \right) \right)}^{\prime }} \right){{b}^{m}}\left( x \right)=\frac{m+n}{m}\left[ {{x}^{n}} \right]{{b}^{m}}\left( x \right)=\left[ {{x}^{n}} \right]{{a}^{m+n}}\left( x \right),$$
 $$\left[ {{x}^{n}} \right]{{b}^{m}}\left( x \right)=\left[ {{x}^{n}} \right]\left( 1-x{{\left( \log a\left( x \right) \right)}^{\prime }} \right){{a}^{m+n}}\left( x \right)=\frac{m}{m+n}\left[ {{x}^{n}} \right]{{a}^{m+n}}\left( x \right).$$
Denote
$${{\left( 1,{{a}^{-\beta }}\left( x \right) \right)}^{-1}}=\left( 1,{}_{\left( \beta  \right)}{{a}^{\beta }}\left( x \right) \right).$$
Then
$$\left[ {{x}^{n}} \right]{}_{\left( \beta  \right)}{{a}^{\beta m}}\left( x \right)=\frac{\beta m}{\beta m+\beta n}\left[ {{x}^{n}} \right]{{a}^{\beta m+\beta n}}\left( x \right).$$
Let ${{{s}_{n}}\left( \varphi  \right)}/{n!}\;$ are the convolution polynomials of the series $a\left( x \right)$. Then
$$_{\left( \beta  \right)}{{a}^{\varphi }}\left( x \right)=\sum\limits_{n=0}^{\infty }{\frac{\varphi }{\varphi +\beta n}}\frac{{{s}_{n}}\left( \varphi +\beta n \right)}{n!}{{x}^{n}}.$$

We note that the identity 
$$\left[ n,\to  \right]{{\left( 1-x{{\left( \log a\left( x \right) \right)}^{\prime }},{{a}^{-1}}\left( x \right) \right)}^{-1}}=\left[ n,\to  \right]\left( {{a}^{n}}\left( x \right),1 \right)$$
is equivalent to the Lagrange expansion formula for an arbitrary series $f\left( x \right)$:
$$\frac{f\left( x \right)}{1-x{{\left( \log a\left( x \right) \right)}^{\prime }}}=\sum\limits_{n=0}^{\infty }{\frac{{{x}^{n}}}{{{a}^{n}}\left( x \right)}}\left[ {{x}^{n}} \right]f\left( x \right){{a}^{n}}\left( x \right).$$
\section{Riordan-Dirichlet matrices}

In this section, we construct the matrices of certain operators in the space of formal Dirichlet series. We will associate the columns of matrices with the generating functions of their elements, i.e. with the formal Dirichlet  series. Numbering of the rows and columns begins with $1$. Coefficient operator will be denoted by $\left[ {{n}^{-s}} \right]$:  $\left[ {{n}^{-s}} \right]a\left( s \right)={{a}_{n}}$. 

Matrix $n$th column of which has the generating function ${{n}^{-s}}a\left( s \right)$ will be  denoted by $\left\langle a\left( s \right),1 \right\rangle $:
$$\left\langle a\left( s \right),1 \right\rangle =\left( \begin{matrix}
   {{a}_{1}} & 0 & 0 & 0 & 0 & 0 & 0 & 0 & \cdots   \\
   {{a}_{2}} & {{a}_{1}} & 0 & 0 & 0 & 0 & 0 & 0 & \cdots   \\
   {{a}_{3}} & 0 & {{a}_{1}} & 0 & 0 & 0 & 0 & 0 & \cdots   \\
   {{a}_{4}} & {{a}_{2}} & 0 & {{a}_{1}} & 0 & 0 & 0 & 0 & \cdots   \\
   {{a}_{5}} & 0 & 0 & 0 & {{a}_{1}} & 0 & 0 & 0 & \cdots   \\
   {{a}_{6}} & {{a}_{3}} & {{a}_{2}} & 0 & 0 & {{a}_{1}} & 0 & 0 & \cdots   \\
   {{a}_{7}} & 0 & 0 & 0 & 0 & 0 & {{a}_{1}} & 0 & \cdots   \\
   {{a}_{8}} & {{a}_{4}} & 0 & {{a}_{2}} & 0 & 0 & 0 & {{a}_{1}} & \cdots   \\
   \vdots  & \vdots  & \vdots  & \vdots  & \vdots  & \vdots  & \vdots  & \vdots  & \ddots   \\
\end{matrix} \right).$$
Then
$$\left\langle a\left( s \right),1 \right\rangle b\left( s \right)=\sum\limits_{n=1}^{\infty }{{{b}_{n}}{{n}^{-s}}a\left( s \right)}=\sum\limits_{n=1}^{\infty }{{{n}^{-s}}\sum\limits_{d|n}^{n}{{{b}_{d}}{{a}_{{n}/{d}\;}}}}=a\left( s \right)b\left( s \right),$$
where symbol $d|n$ means that the summation is over all natural divisors $d$ of the number $n$. Thus, matrix $\left\langle a\left( s \right),1 \right\rangle $ corresponds to the operator of multiplication by the  series $a\left( s \right)$. Set of all such matrices form an algebra isomorphic to the algebra of formal Dirichlet series:
$$\left\langle a\left( s \right),1 \right\rangle +\left\langle b\left( s \right),1 \right\rangle =\left\langle a\left( s \right)+b\left( s \right),1 \right\rangle ,$$
$$\left\langle a\left( s \right),1 \right\rangle \left\langle b\left( s \right),1 \right\rangle =\left\langle b\left( s \right),1 \right\rangle \left\langle a\left( s \right),1 \right\rangle =\left\langle a\left( s \right)b\left( s \right),1 \right\rangle .$$
{\bfseries Remark 3.}
Matrices $\left\langle a\left( s \right),1 \right\rangle $, denoted by $A=\left[ \left( {{a}_{n}} \right) \right]$ , were introduced in [8], [9]. In [9], [10], they are considered as matrices of  the operators in the Hilbert space and are called (in [10]) the Dirichlet matrices, or $D$-matrices.

In the algebra of formal Dirichlet series, the power and logarithm of the series $a\left( s \right)$, ${{a}_{1}}=1$, are defined as
$${{a}^{\varphi }}\left( s \right)=\sum\limits_{n=0}^{\infty }{\left( \begin{matrix}
   \varphi   \\
   n  \\
\end{matrix} \right){{\left( a\left( s \right)-1 \right)}^{n}}}, \qquad\log a\left( s \right)=\sum\limits_{n=1}^{\infty }{\frac{{{\left( -1 \right)}^{n-1}}}{n}}{{\left( a\left( s \right)-1 \right)}^{n}}.$$
Power of the series can also be defined as
$${{a}^{\varphi }}\left( s \right)=\exp \left( \varphi \log a\left( s \right) \right)=\sum\limits_{n=0}^{\infty }{\frac{{{\varphi }^{n}}}{n!}}{{\left( \log a\left( s \right) \right)}^{n}}=\sum\limits_{n=1}^{\infty }{{{h}_{n}}\left( \varphi  \right)}{{n}^{-s}},$$
where ${{h}_{n}}\left( \varphi  \right)$ are the polynomials in $\varphi $ of degree  $<n$ which, like in the case of power series, can be called convolution polynomials. Explicit form of these polynomials:
$${{h}_{1}}\left( \varphi  \right)=1,  \qquad{{h}_{n}}\left( \varphi  \right)=\sum\limits_{m=1}^{n}{{{\varphi }^{m}}\sum{\frac{b_{2}^{{{m}_{2}}}b_{3}^{{{m}_{3}}}...\text{ }b_{n}^{{{m}_{n}}}}{{{m}_{2}}!{{m}_{3}}!\text{ }...\text{ }{{m}_{n}}!}}}\text{ },$$
where ${{b}_{i}}=\left[ {{x}^{i}} \right]\log a\left( s \right)$ and summation of the coefficient of ${{\varphi }^{m}}$ is over all monomials $b_{2}^{{{m}_{2}}}b_{3}^{{{m}_{3}}}...\text{ }b_{n}^{{{m}_{n}}}$ for which $\prod\nolimits_{i=2}^{n}{{{i}^{{{m}_{i}}}}}=n$, $\sum\nolimits_{i=2}^{n}{{{m}_{i}}}=m$.

Matrix $n$th column of which has the generating function ${{n}^{-s}}{{a}^{\ln n}}\left( s \right)$ will be  denoted by $\left\langle 1,a\left( s \right) \right\rangle $:
$$\left\langle 1,a\left( s \right) \right\rangle =\left( \begin{matrix}
   1 & 0 & 0 & 0 & 0 & 0 & 0 & 0 & \cdots   \\
   0 & a_{1}^{\left( 2 \right)} & 0 & 0 & 0 & 0 & 0 & 0 & \cdots   \\
   0 & 0 & a_{1}^{\left( 3 \right)} & 0 & 0 & 0 & 0 & 0 & \cdots   \\
   0 & a_{2}^{\left( 2 \right)} & 0 & a_{1}^{\left( 4 \right)} & 0 & 0 & 0 & 0 & \cdots   \\
   0 & 0 & 0 & 0 & a_{1}^{\left( 5 \right)} & 0 & 0 & 0 & \cdots   \\
   0 & a_{3}^{\left( 2 \right)} & a_{2}^{\left( 3 \right)} & 0 & 0 & a_{1}^{\left( 6 \right)} & 0 & 0 & \cdots   \\
   0 & 0 & 0 & 0 & 0 & 0 & a_{1}^{\left( 7 \right)} & 0 & \cdots   \\
   0 & a_{4}^{\left( 2 \right)} & 0 & a_{2}^{\left( 4 \right)} & 0 & 0 & 0 & a_{1}^{\left( 8 \right)} & \cdots   \\
   \vdots  & \vdots  & \vdots  & \vdots  & \vdots  & \vdots  & \vdots  & \vdots  & \ddots   \\
\end{matrix} \right),$$
where
$$\left[ {{n}^{-s}} \right]{{a}^{\ln m}}\left( s \right)=a_{n}^{\left( m \right)},  \qquad a_{1}^{\left( 1 \right)}=1,  \qquad a_{n}^{\left( 1 \right)}=0.$$
Then
$$\left\langle 1,a\left( s \right) \right\rangle b\left( s \right)=\sum\limits_{n=1}^{\infty }{{{b}_{n}}{{n}^{-s}}{{a}^{\ln n}}\left( s \right)=}\sum\limits_{n=1}^{\infty }{{{n}^{-s}}\sum\limits_{d|n}^{n}{{{b}_{d}}a_{{n}/{d}\;}^{\left( d \right)}}}=b\left( s \right)\circ a\left( s \right),$$
where sign $\circ $ means the operation under consideration similar to a composition of power series. Denote
$$\left\langle b\left( s \right),1 \right\rangle \left\langle 1,a\left( s \right) \right\rangle =\left\langle b\left( s \right),a\left( s \right) \right\rangle .$$
 {\bfseries Theorem 2.} \emph{Matrices $\left\langle b\left( s \right),a\left( s \right) \right\rangle $, ${{b}_{1}}\ne 0$, ${{a}_{1}}\ne 0$, form a group with respect to the matrix multiplication whose elements are multiplied by the rule:}
$$\left\langle b\left( s \right),a\left( s \right) \right\rangle \left\langle f\left( s \right),g\left( s \right) \right\rangle =\left\langle b\left( s \right)\left( f\left( s \right)\circ a\left( s \right) \right),a\left( x \right)\left( g\left( s \right)\circ a\left( s \right) \right) \right\rangle .$$
{\bfseries Proof.}
Since
$$\left\langle 1,a\left( s \right) \right\rangle {{m}^{-s}}b\left( s \right)=\sum\limits_{n=1}^{\infty }{{{b}_{n}}{{\left( mn \right)}^{-s}}{{a}^{\ln mn}}\left( s \right)}={{m}^{-s}}{{a}^{\ln m}}\left( s \right)\sum\limits_{n=1}^{\infty }{{{b}_{n}}}{{n}^{-s}}{{a}^{\ln n}}\left( s \right),$$
then 
$$\left\langle 1,a\left( s \right) \right\rangle \left\langle b\left( s \right),1 \right\rangle =\left\langle b\left( s \right)\circ a\left( s \right),a\left( s \right) \right\rangle .$$
Hence,
$$\left\langle 1,a\left( s \right) \right\rangle b\left( s \right)c\left( s \right)=\left( b\left( s \right)\circ a\left( s \right) \right)\left( c\left( s \right)\circ a\left( s \right) \right),$$
$$\left\langle 1,a\left( s \right) \right\rangle {{b}^{n}}\left( s \right)={{\left( b\left( s \right)\circ a\left( s \right) \right)}^{n}},$$
and by definition of power of series
$$\left\langle 1,a\left( s \right) \right\rangle {{b}^{\varphi }}\left( s \right)={{\left( b\left( s \right)\circ a\left( s \right) \right)}^{\varphi }}.$$
Then
$$\left\langle 1,a\left( s \right) \right\rangle {{m}^{-s}}{{b}^{\ln m}}\left( s \right)={{m}^{-s}}{{a}^{\ln m}}\left( s \right){{\left( b\left( x \right)\circ a\left( x \right) \right)}^{\ln m}},$$
$$\left\langle 1,a\left( s \right) \right\rangle \left\langle 1,b\left( s \right) \right\rangle =\left\langle 1,a\left( s \right)\left( b\left( s \right)\circ a\left( s \right) \right) \right\rangle .$$
If
$${{\left\langle 1,a\left( s \right) \right\rangle }^{-1}}{{b}^{-1}}\left( s \right)=f\left( s \right),   \qquad{{\left\langle 1,a\left( s \right) \right\rangle }^{-1}}{{a}^{-1}}\left( s \right)=g\left( s \right),$$
then
$${{\left\langle b\left( s \right),a\left( s \right) \right\rangle }^{-1}}=\left\langle f\left( s \right),g\left( s \right) \right\rangle .$$

Matrix of the differential operator in the space of formal Dirichlet series will be denoted by ${{D}_{s}}$:
$${{D}_{s}}a\left( s \right)={a}'\left( s \right)=\sum\limits_{n=1}^{\infty }{\ln \left( {1}/{n}\; \right){{a}_{n}}}{{n}^{-s}}.$$
Then
$${{\left( a\left( s \right)b\left( s \right) \right)}^{\prime }}=a\left( s \right){b}'\left( s \right)+{a}'\left( s \right)b\left( s \right),  \qquad{{\left( {{a}^{n}}\left( s \right) \right)}^{\prime }}=n{{a}^{n-1}}\left( s \right){a}'\left( s \right),$$
$${{\left( {{a}^{\varphi }}\left( s \right) \right)}^{\prime }}=\varphi {a}'\left( s \right)\sum\limits_{n=1}^{\infty }{\left( \begin{matrix}
   \varphi -1  \\
   n-1  \\
\end{matrix} \right){{\left( a\left( s \right)-1 \right)}^{n-1}}}=\varphi {{a}^{\varphi -1}}\left( s \right){a}'\left( s \right),$$
$${{\left( \log a\left( s \right) \right)}^{\prime }}={a}'\left( s \right)\sum\limits_{n=1}^{\infty }{{{\left( -1 \right)}^{n-1}}{{\left( a\left( s \right)-1 \right)}^{n-1}}}={a}'\left( s \right){{a}^{-1}}\left( s \right).$$
{\bfseries Theorem 3.} \emph{Each formal Dirichlet series $a\left( s \right)$, ${{a}_{1}}=1$, is associated  by means of  the transform
$${{a}^{\varphi }}\left( s \right)=\sum\limits_{n=1}^{\infty }{\frac{{{n}^{-s}}}{{{a}^{\beta \ln n}}\left( s \right)}\left[ {{n}^{-s}} \right]}\left( 1+\beta {{\left( \log a\left( s \right) \right)}^{\prime }} \right){{a}^{\varphi +\beta \ln n}}\left( s \right)$$
with the set of series $_{\left( \beta  \right)}a\left( s \right)$, $_{\left( 0 \right)}a\left( s \right)=a\left( s \right)$, such that}
$$\left[ {{n}^{-s}} \right]{}_{\left( \beta  \right)}{{a}^{\varphi }}\left( s \right)=\left[ {{n}^{-s}} \right]\left( 1+\beta {{\left( \log a\left( s \right) \right)}^{\prime }} \right){{a}^{\varphi +\beta \ln n}}\left( s \right)=\frac{\varphi }{\varphi +\beta \ln n}\left[ {{n}^{-s}} \right]{{a}^{\varphi +\beta \ln n}}\left( s \right),$$
$$\left[ {{n}^{-s}} \right]\left( 1-\beta {{\left( \log {}_{\left( \beta  \right)}a\left( s \right) \right)}^{\prime }} \right){}_{\left( \beta  \right)}{{a}^{\varphi }}\left( s \right)=\frac{\varphi +\beta \ln n}{\varphi }\left[ {{n}^{-s}} \right]{}_{\left( \beta  \right)}{{a}^{\varphi }}\left( s \right)=\left[ {{n}^{-s}} \right]{{a}^{\varphi +\beta \ln n}}\left( s \right).$$
{\bfseries Proof.}
If the matrices $\left\langle 1,{{a}^{-1}}\left( s \right) \right\rangle $, ${{a}_{1}}=1$, $\left\langle 1,b\left( s \right) \right\rangle $, ${{b}_{1}}=1$, are mutually inverse, then
$$\left\langle 1,{{a}^{-1}}\left( s \right) \right\rangle b\left( s \right)=a\left( s \right),   \qquad\left\langle 1,b\left( s \right) \right\rangle a\left( s \right)=b\left( s \right).$$
Since
$${{\left( {{n}^{-s}}{{b}^{\ln n}}\left( s \right) \right)}^{\prime }}={{n}^{-s}}{{b}^{\ln n}}\left( s \right)\ln \left( {1}/{n}\; \right)\left( 1-{b}'\left( s \right){{b}^{-1}}\left( s \right) \right),$$
then
$${{D}_{s}}\left\langle 1,b\left( s \right) \right\rangle =\left\langle 1-{{\left( \log b\left( s \right) \right)}^{\prime }},b\left( s \right) \right\rangle {{D}_{s}},$$
$$\left\langle 1,b\left( s \right) \right\rangle {a}'\left( s \right)=\frac{{b}'\left( s \right)}{1-{{\left( \log b\left( s \right) \right)}^{\prime }}}.$$
From this we find:
$${{\left\langle 1-{{\left( \log b\left( s \right) \right)}^{\prime }},b\left( s \right) \right\rangle }^{-1}}=\left\langle 1+{{\left( \log a\left( s \right) \right)}^{\prime }},{{a}^{-1}}\left( s \right) \right\rangle .$$
Denote
$$\left[ {{n}^{-s}} \right]{{a}^{\ln m}}\left( s \right)=a_{n}^{\left( m \right)},    \qquad\left[ {{n}^{-s}} \right]\left( 1+{{\left( \log a\left( s \right) \right)}^{\prime }} \right){{a}^{\ln m}}\left( s \right)=c_{n}^{\left( m \right)},$$
$${{a}_{m}}\left( s \right)=\sum\limits_{n=1}^{\infty }{a_{n}^{\left( mn \right)}}{{n}^{-s}},   \qquad{{c}_{m}}\left( s \right)=\sum\limits_{n=1}^{\infty }{c_{n}^{\left( mn \right)}}{{n}^{-s}}.$$
We construct the matrix $A$ $m$th column of which has the generating function ${{m}^{-s}}{{a}_{m}}\left( s \right)$ and the matrix $C$ $m$th column of which has the generating function ${{m}^{-s}}{{c}_{m}}\left( s \right)$:
$$\left( \begin{matrix}
   a_{1}^{\left( 1 \right)} & 0 & 0 & 0 & 0 & 0 & \cdots   \\
   a_{2}^{\left( 2 \right)} & a_{1}^{\left( 2 \right)} & 0 & 0 & 0 & 0 & \cdots   \\
   a_{3}^{\left( 3 \right)} & 0 & a_{1}^{\left( 3 \right)} & 0 & 0 & 0 & \cdots   \\
   a_{4}^{\left( 4 \right)} & a_{2}^{\left( 4 \right)} & 0 & a_{1}^{\left( 4 \right)} & 0 & 0 & \cdots   \\
   a_{5}^{\left( 5 \right)} & 0 & 0 & 0 & a_{1}^{\left( 5 \right)} & 0 & \cdots   \\
   a_{6}^{\left( 6 \right)} & a_{3}^{\left( 6 \right)} & a_{2}^{\left( 6 \right)} & 0 & 0 & a_{1}^{\left( 6 \right)} & \cdots   \\
   \vdots  & \vdots  & \vdots  & \vdots  & \vdots  & \vdots  & \ddots   \\
\end{matrix} \right),\quad\left( \begin{matrix}
   c_{1}^{\left( 1 \right)} & 0 & 0 & 0 & 0 & 0 & \cdots   \\
   c_{2}^{\left( 2 \right)} & c_{1}^{\left( 2 \right)} & 0 & 0 & 0 & 0 & \cdots   \\
   c_{3}^{\left( 3 \right)} & 0 & c_{1}^{\left( 3 \right)} & 0 & 0 & 0 & \cdots   \\
   c_{4}^{\left( 4 \right)} & c_{2}^{\left( 4 \right)} & 0 & c_{1}^{\left( 4 \right)} & 0 & 0 & \cdots   \\
   c_{5}^{\left( 5 \right)} & 0 & 0 & 0 & c_{1}^{\left( 5 \right)} & 0 & \cdots   \\
   c_{6}^{\left( 6 \right)} & c_{3}^{\left( 6 \right)} & c_{2}^{\left( 6 \right)} & 0 & 0 & c_{1}^{\left( 6 \right)} & \cdots   \\
   \vdots  & \vdots  & \vdots  & \vdots  & \vdots  & \vdots  & \ddots   \\
\end{matrix} \right)$$
It's obvious that
$$\left[ n,\to  \right]A=\left[ n,\to  \right]\left\langle {{a}^{\ln n}}\left( s \right),1 \right\rangle ,$$
$$\left[ n,\to  \right]C=\left[ n,\to  \right]\left\langle 1+{{\left( \log a\left( s \right) \right)}^{\prime }}{{a}^{\ln n}}\left( s \right),1 \right\rangle .$$
Since
$$\left( 1+{a}'\left( s \right){{a}^{-1}}\left( s \right) \right){{a}^{\ln m}}\left( s \right)={{a}^{\ln m}}\left( s \right)+\frac{1}{\ln m}{{\left( {{a}^{\ln m}}\left( s \right) \right)}^{\prime }},$$
or
$$\left[ {{n}^{-s}} \right]\left( 1+{{\left( \log a\left( s \right) \right)}^{\prime }} \right){{a}^{\ln m}}\left( s \right)=\frac{\ln \left( {m}/{n}\; \right)}{\ln m}\left[ {{n}^{-s}} \right]{{a}^{\ln m}}\left( s \right),$$
then
$$\left[ {{\left( nm \right)}^{-s}} \right]A{{m}^{-s}}\left( 1+{{\left( \log a\left( s \right) \right)}^{\prime }} \right){{a}^{-\ln m}}\left( s \right)=\left[ {{\left( nm \right)}^{-s}} \right]C{{m}^{-s}}{{a}^{-\ln m}}\left( s \right)=$$
$$=\left[ {{n}^{-s}} \right]\left( 1+{{\left( \log a\left( x \right) \right)}^{\prime }} \right){{a}^{\ln n}}\left( s \right)=1, n=1; \qquad=0, n>1.$$
Thus, 
$$A=\left\langle 1-{{\left( \log b\left( s \right) \right)}^{\prime }},b\left( s \right) \right\rangle ,  \qquad C=\left\langle 1,b\left( s \right) \right\rangle ,$$
$$\left[ {{n}^{-s}} \right]\left( 1-{{\left( \log b\left( s \right) \right)}^{\prime }} \right){{b}^{\ln m}}\left( s \right)=\frac{\ln mn}{\ln m}\left[ {{n}^{-s}} \right]{{b}^{\ln m}}\left( s \right)=\left[ {{n}^{-s}} \right]{{a}^{\ln mn}}\left( s \right),$$
$$\left[ {{n}^{-s}} \right]{{b}^{\ln m}}\left( s \right)=\left[ {{n}^{-s}} \right]\left( 1+{{\left( \log a\left( s \right) \right)}^{\prime }} \right){{a}^{\ln mn}}\left( s \right)=\frac{\ln m}{\ln mn}\left[ {{n}^{-s}} \right]{{a}^{\ln mn}}\left( s \right).$$
Denote
$${{\left\langle 1,{{a}^{-\beta }}\left( s \right) \right\rangle }^{-1}}=\left\langle 1,{}_{\left( \beta  \right)}{{a}^{\beta }}\left( s \right) \right\rangle.$$
Then
$$\left[ {{n}^{-s}} \right]{}_{\left( \beta  \right)}{{a}^{\beta \ln m}}\left( s \right)=\frac{\beta \ln m}{\beta \ln m+\beta \ln n}\left[ {{n}^{-s}} \right]{{a}^{\beta \ln m+\beta \ln n}}\left( s \right).$$
Let ${{h}_{n}}\left( \varphi  \right)$ are the convolution polynomials of the series $a\left( s \right)$. Then
$$_{\left( \beta  \right)}{{a}^{\varphi }}\left( s \right)=\sum\limits_{n=1}^{\infty }{\frac{\varphi }{\varphi +\beta \ln n}}{{h}_{n}}\left( \varphi +\beta \ln n \right){{n}^{-s}}.$$

We note that the identity 
$$\left[ n,\to  \right]{{\left\langle 1+{{\left( \log a\left( s \right) \right)}^{\prime }},{{a}^{-1}}\left( s \right) \right\rangle }^{-1}}=\left[ n,\to  \right]\left\langle {{a}^{\ln n}}\left( s \right),1 \right\rangle $$
is equivalent to a formula similar to the Lagrange expansion formula:
$$\frac{f\left( s \right)}{1+{{\left( \log a\left( s \right) \right)}^{\prime }}}=\sum\limits_{n=1}^{\infty }{\frac{{{n}^{-s}}}{{{a}^{\ln n}}\left( s \right)}}\left[ {{n}^{-s}} \right]f\left( s \right){{a}^{\ln n}}\left( s \right).$$
\section{Some examples}
In this section, we show how matrices $\left\langle 1,a\left( s \right) \right\rangle $ can be used for the next construction. Each formal power series $a\left( x \right)$, ${{a}_{0}}=1$, with convolution polynomials ${{{s}_{n}}\left( \varphi  \right)}/{n!}\;$  is put in correspondence with the formal Dirichlet series $a\left( s \right)$  such that
$${{a}^{\varphi }}\left( s \right)=\prod\limits_{p=2}^{\infty }{\left( \sum\limits_{n=0}^{\infty }{\frac{{{s}_{n}}\left( \varphi  \right)}{n!}{{p}^{-ns}}} \right)}=1+\sum\limits_{n=2}^{\infty }{\frac{{{s}_{{{m}_{1}}}}\left( \varphi  \right){{s}_{{{m}_{2}}}}\left( \varphi  \right)...{{s}_{{{m}_{r}}}}\left( \varphi  \right)}{{{m}_{1}}!{{m}_{2}}!\text{ }...\text{ }{{m}_{r}}!}}{{n}^{-s}},$$  
where the product is taken over all prime numbers and $n=p_{1}^{{{m}_{1}}}p_{2}^{{{m}_{2}}}...\text{ }p_{r}^{{{m}_{r}}}$ is the canonical decomposition of  number $n$. With this correspondence, the group of series $a\left( x \right)$ is isomorphic to the group of  series $a\left( s \right)$: if $a\left( x \right)b\left( x \right)=c\left( x \right)$, then $a\left( s \right)b\left( s \right)=c\left( s \right)$. If $\left[ {{x}^{m}} \right]\log a\left( x \right)={{b}_{m}}$, then $\left[ {{n}^{-s}} \right]\log a\left( s \right)={{b}_{m}}$ when $n={{p}^{m}}$, where $p$ is the  prime number, and $\left[ {{n}^{-s}} \right]\log a\left( s \right)=0$ otherwise. Polynomials
$${{h}_{n}}\left( \varphi  \right)=\frac{{{s}_{{{m}_{1}}}}\left( \varphi  \right){{s}_{{{m}_{2}}}}\left( \varphi  \right)...{{s}_{{{m}_{r}}}}\left( \varphi  \right)}{{{m}_{1}}!{{m}_{2}}!\text{ }...\text{ }{{m}_{r}}!}$$
are the convolution polynomials of the series $a\left( s \right)$. Dirichlet series corresponding to the exponential series will be denoted by $\varepsilon \left( s \right)$:
$${{\varepsilon }^{\varphi }}\left( s \right)=1+\sum\limits_{n=2}^{\infty }{\frac{{{\varphi }^{{{m}_{1}}+{{m}_{2}}+...+{{m}_{r}}}}}{{{m}_{1}}!{{m}_{2}}!\text{ }...\text{ }{{m}_{r}}!}{{n}^{-s}}}.$$
Note that $\log \varepsilon \left( s \right)=\sum{{{p}^{-s}}}$, where the summation is over all prime numbers. Let $_{\left( 1 \right)}\varepsilon \left( s \right)$ denote the series associated with $\varepsilon \left( s \right)$ by Theorem 3. Then
$$\left[ {{n}^{-s}} \right]{{\varepsilon }^{\varphi }}\left( s \right)=\frac{{{\varphi }^{v\left( n \right)}}}{f\left( n \right)},  \qquad\left[ {{n}^{-s}} \right]{}_{\left( 1 \right)}{{\varepsilon }^{\varphi }}\left( s \right)=\frac{\varphi {{\left( \varphi +\ln n \right)}^{v\left( n \right)-1}}}{f\left( n \right)},$$
where 
$$v\left( 1 \right)=0,  \quad v\left( n \right)={{m}_{1}}+{{m}_{2}}+...+{{m}_{r}},  \qquad f\left( 1 \right)=1,  \quad f\left( n \right)={{m}_{1}}!{{m}_{2}}!...{{m}_{r}}!.$$
From
$$_{\left( 1 \right)}{{\varepsilon }^{\varphi +\beta }}\left( s \right)={}_{\left( 1 \right)}{{\varepsilon }^{\varphi }}\left( s \right){}_{\left( 1 \right)}{{\varepsilon }^{\beta }}\left( s \right)$$
we obtain an analog of the Abel's generalized binomial formula:
$$\left( \varphi +\beta  \right){{\left( \varphi +\beta +\ln n \right)}^{v\left( n \right)-1}}=\sum\limits_{d|n}{{{\left( \begin{matrix}
   n  \\
   d  \\
\end{matrix} \right)}_{f}}}\varphi {{\left( \varphi +\ln d \right)}^{v\left( d \right)-1}}\beta {{\left( \beta +\ln \left( {n}/{d}\; \right) \right)}^{v\left( {n}/{d}\; \right)-1}},$$
where
$${{\left( \begin{matrix}
   n  \\
   d  \\
\end{matrix} \right)}_{f}}=\frac{f\left( n \right)}{f\left( d \right)f\left( {n}/{d}\; \right)}.$$ 
Since 
$$\left[ {{n}^{-s}} \right]\left( 1-{{\left( \log {}_{\left( 1 \right)}\varepsilon \left( s \right) \right)}^{\prime }} \right){}_{\left( 1 \right)}{{\varepsilon }^{\varphi }}\left( s \right)=\left[ {{n}^{-s}} \right]{{\varepsilon }^{\varphi +\ln n}}\left( s \right),$$
from
$$\left( 1-{{\left( \log {}_{\left( 1 \right)}\varepsilon \left( s \right) \right)}^{\prime }} \right){}_{\left( 1 \right)}{{\varepsilon }^{\varphi +\beta }}\left( s \right)=\left( 1-{{\left( \log {}_{\left( 1 \right)}\varepsilon \left( s \right) \right)}^{\prime }} \right){}_{\left( 1 \right)}{{\varepsilon }^{\varphi }}\left( s \right){}_{\left( 1 \right)}{{\varepsilon }^{\beta }}\left( s \right)$$
we obtain 
$${{\left( \varphi +\beta +\ln n \right)}^{v\left( n \right)}}=\sum\limits_{d|n}{{{\left( \begin{matrix}
   n  \\
   d  \\
\end{matrix} \right)}_{f}}}{{\left( \varphi +\ln d \right)}^{v\left( d \right)}}\beta {{\left( \beta +\ln \left( {n}/{d}\; \right) \right)}^{v\left( {n}/{d}\; \right)-1}}.$$
Since
$$\left[ {{n}^{-s}} \right]\left\langle 1,a\left( s \right) \right\rangle b\left( s \right)=\sum\limits_{d|n}^{n}{{{b}_{d}}a_{{n}/{d}\;}^{\left( d \right)}},  \qquad a_{{n}/{d}\;}^{\left( d \right)}=\left[ {{\left( {n}/{d}\; \right)}^{-s}} \right]{{a}^{\ln d}}\left( s \right),$$
from
$$_{\left( 1 \right)}{{\varepsilon }^{\varphi }}\left( s \right)=\left\langle 1,{}_{\left( 1 \right)}\varepsilon \left( s \right) \right\rangle {{\varepsilon }^{\varphi }}\left( s \right),  \qquad {{\varepsilon }^{\varphi }}\left( x \right)=\left\langle 1,{{\varepsilon }^{-1}}\left( s \right) \right\rangle {}_{\left( 1 \right)}{{\varepsilon }^{\varphi }}\left( s \right)$$
we obtain 
$$\varphi {{\left( \varphi +\ln n \right)}^{v\left( n \right)-1}}=\sum\limits_{d|n}{{{\left( \begin{matrix}
   n  \\
   d  \\
\end{matrix} \right)}_{f}}{{\varphi }^{v\left( d \right)}}\ln d{{\left( \ln n \right)}^{v\left( {n}/{d}\; \right)-1}}},$$
$${{\varphi }^{v\left( n \right)}}=\sum\limits_{d|n}{{{\left( \begin{matrix}
   n  \\
   d  \\
\end{matrix} \right)}_{f}}\varphi {{\left( \varphi +\ln d \right)}^{v\left( d \right)-1}}{{\left( \ln \left( {1}/{d}\; \right) \right)}^{v\left( {n}/{d}\; \right)}}}.$$
When $n={{p}^{m}}$, where $p$ is the  prime number, obtained formulas take the form of Abel identities [3],  [11, pp. 92-99]:
$$\left( \varphi +\beta  \right){{\left( \varphi +\beta +ma \right)}^{m-1}}=\sum\limits_{k=0}^{m}{\left( \begin{matrix}
   m  \\
   k  \\
\end{matrix} \right)\varphi {{\left( \varphi +ka \right)}^{k-1}}\beta {{\left( \beta +\left( m-k \right)a \right)}^{m-k-1}}},$$
$${{\left( \varphi +\beta +ma \right)}^{m}}=\sum\limits_{k=0}^{m}{\left( \begin{matrix}
   m  \\
   k  \\
\end{matrix} \right){{\left( \varphi +ka \right)}^{k}}}\beta {{\left( \beta +\left( m-k \right)a \right)}^{m-k-1}},$$
$$\varphi {{\left( \varphi +ma \right)}^{m-1}}=\sum\limits_{k=0}^{m}{\left( \begin{matrix}
   m  \\
   k  \\
\end{matrix} \right)}{{\varphi }^{k}}ka{{\left( ma \right)}^{m-k-1}},$$
$${{\varphi }^{m}}=\sum\limits_{k=0}^{m}{\left( \begin{matrix}
   m  \\
   k  \\
\end{matrix} \right)\varphi {{\left( \varphi +ka \right)}^{k-1}}}{{\left( -ka \right)}^{m-k}},  \qquad a=\ln p.$$

We generalize this example. Let there be given the series $a\left( s \right)$ such that
$$\left[ {{n}^{-s}} \right]{{a}^{\varphi }}\left( s \right)=\frac{{{u}_{n}}\left( \varphi  \right)}{f\left( n \right)},$$
where
$${{u}_{n}}\left( \varphi  \right)={{s}_{{{m}_{1}}}}\left( \varphi  \right){{s}_{{{m}_{2}}}}\left( \varphi  \right)...{{s}_{{{m}_{r}}}}\left( \varphi  \right),  \qquad n=p_{1}^{{{m}_{1}}}p_{2}^{{{m}_{2}}}...\text{ }p_{r}^{{{m}_{r}}}.$$
Then
$$\left[ {{n}^{-s}} \right]{}_{\left( \beta  \right)}{{a}^{\varphi }}\left( s \right)=\frac{\varphi }{\varphi +\beta \ln n}\frac{{{u}_{n}}\left( \varphi +\beta \ln n \right)}{f\left( n \right)},$$
$$_{\left( \beta  \right)}{{a}^{\varphi }}\left( s \right)=\left\langle 1,{}_{\left( \beta  \right)}{{a}^{\beta }}\left( s \right) \right\rangle {{a}^{\varphi }}\left( s \right),  \qquad{{a}^{\varphi }}\left( s \right)=\left\langle 1,{{a}^{-\beta }}\left( s \right) \right\rangle {}_{\left( \beta  \right)}{{a}^{\varphi }}\left( s \right),$$
$$\frac{\varphi }{\varphi +\beta \ln n}{{u}_{n}}\left( \varphi +\beta \ln n \right)=\sum\limits_{d|n}{{{\left( \begin{matrix}
   n  \\
   d  \\
\end{matrix} \right)}_{f}}{{u}_{d}}\left( \varphi  \right)\frac{\ln d}{\ln n}}{{u}_{{n}/{d}\;}}\left( \beta \ln n \right),$$
$${{u}_{n}}\left( \varphi  \right)=\sum\limits_{d|n}{{{\left( \begin{matrix}
   n  \\
   d  \\
\end{matrix} \right)}_{f}}\frac{\varphi }{\varphi +\beta \ln d}{{u}_{d}}\left( \varphi +\beta \ln d \right)}{{u}_{{n}/{d}\;}}\left( \beta \ln \left( {1}/{d}\; \right) \right).$$
Since ${{u}_{{{p}^{m}}}}\left( \varphi  \right)={{s}_{m}}\left( \varphi  \right)$, when $n={{p}^{m}}$ formulas take the form of mutually inverse relations for the Lagrange series:
$$\frac{\varphi }{\varphi +ma}{{s}_{m}}\left( \varphi +ma \right)=\sum\limits_{k=0}^{m}{\left( \begin{matrix}
   m  \\
   k  \\
\end{matrix} \right)}{{s}_{k}}\left( \varphi  \right)\frac{k}{m}{{s}_{m-k}}\left( ma \right),$$
$${{s}_{m}}\left( \varphi  \right)=\sum\limits_{k=0}^{m}{\left( \begin{matrix}
   m  \\
   k  \\
\end{matrix} \right)}\frac{\varphi }{\varphi +ka}{{s}_{k}}\left( \varphi +ka \right){{s}_{m-k}}\left( -ka \right),  \qquad a=\beta \ln p.$$

We note the identities for the coefficients ${n\choose d}_f$ similar to the identities  
$$\sum\limits_{k=0}^{n}{\left( \begin{matrix}
   n  \\
   k  \\
\end{matrix} \right)}={{2}^{n}};  \qquad\sum\limits_{k=0}^{n}{\left( \begin{matrix}
   n  \\
   k  \\
\end{matrix} \right)k{{\left( -1 \right)}^{n-k}}}=0, \quad n\ne 1.$$
Since $\varepsilon \left( s \right)\varepsilon \left( s \right)={{\varepsilon }^{2}}\left( s \right)$;  ${\varepsilon }'\left( x \right){{\varepsilon }^{-1}}\left( x \right)={{\left( \log \varepsilon \left( x \right) \right)}^{\prime }}$, then
$$\sum\limits_{d|n}{{{\left( \begin{matrix}
   n  \\
   d  \\
\end{matrix} \right)}_{f}}}={{2}^{v\left( n \right)}};   \qquad\sum\limits_{d|n}{{{\left( \begin{matrix}
   n  \\
   d  \\
\end{matrix} \right)}_{f}}}\ln d{{\left( -1 \right)}^{v\left( {n}/{d}\; \right)}}=0, \quad n\ne p.$$

E-mail: {evgeniy\symbol{"5F}burlachenko@list.ru}
\end{document}